\documentclass[12pt]{amsart}

\usepackage{amsmath}
\usepackage{amsfonts}   
\usepackage{amssymb}

\newcommand{\CC}{\mathbb{C}}

\newcommand{\OO}{\mathcal{O}}

\newcommand{\A}{\mathcal{A}}
\newcommand{\NN}{\mathbf{N}}

\renewcommand{\bar}{\overline}
\newcommand{\Kbar}{{\bar K}}
\renewcommand{\L}{{\mathcal L}}

\newcommand{\hhat}{\hat{h}}

\DeclareMathOperator{\Gal}{Gal}

\newcounter{alphenum}

\newtheorem{thm}{Theorem}[section]

\theoremstyle{definition}

\newtheorem*{defn}{Definition}

\newenvironment{notation}[0]{%
  \begin{list}%
    {}%
    {\setlength{\itemindent}{0pt}
     \setlength{\labelwidth}{4\parindent}
     \setlength{\labelsep}{\parindent}
     \setlength{\leftmargin}{5\parindent}
     \setlength{\itemsep}{0pt}
     }%
   }%
  {\end{list}}

\begin{document}

%%%%%%%%%%%%%%%%%%%%%%%%%%%%%%%%%%%%%%%%%
%               Opening                 %
%%%%%%%%%%%%%%%%%%%%%%%%%%%%%%%%%%%%%%%%%

\title{Equidistribution of small subvarieties of an abelian variety}
\author{Matthew Baker and Su-ion Ih}
  \address{Department of Mathematics\\
        University of Georgia\\
        Athens, GA 30602--7403}
  \email{mbaker@math.uga.edu \\ ih@math.uga.edu}

\thanks{The first author's research was supported in part
by NSF Research Grant DMS-0300784.}

\begin{abstract}
We prove an equidistribution result for small subvarieties of an abelian variety
which generalizes the Szpiro-Ullmo-Zhang theorem on equidistribution of small points.
% It says that small varieties are equidistributed on an abelian variety. 
\end{abstract}

%%%%%%%%%%%%%%%%%%%%%%%%%%%%%%%%%%%%%%%%%
%               Document Text           %
%%%%%%%%%%%%%%%%%%%%%%%%%%%%%%%%%%%%%%%%%

\maketitle

\section{Introduction} \label{sec:intro}

\subsection{Notation}
\label{Notation}

The following notation and conventions will be used throughout this paper:

\begin{notation}
\item[$K$]
a number field.
\item[$\OO_K$]
the ring of integers of $K$.
\item[$A$] 
an abelian variety defined over $K$.

% \item[$X$]
% a projective variety defined over $K$
% \item[$\mathcal X$]
% a model for $X$ over $\OO_K$, i.e., 
% an integral scheme projective and flat over Spec $\OO_K$ whose generic
% fiber is $X$
% \item[$L$]
% a line bundle on $X$
% \item[$\L$]
% a line bundle on $\X$ with generic fiber $L$
% \item[$\overline {L}$]
% a metrized line bundle on $X$, i.e., a line bundle $L$ on $X$,
% together with a smooth hermitian metric on $L_\CC$
% \item[$\overline {\mathcal L}$]
% a metrized line bundle on $\X$, i.e., a line bundle $\L$ on $\X$,
% together with a smooth hermitian metric on $L_\CC$
% \item[$c_1(L)$]
% the first Chern class of $L$, which when $L$ is ample coincides with the
% positive integer $\deg_L(X)$
% \item[$c_1(\overline{L})$]
% the curvature form of $\overline {L}$, which is a
% smooth $(1,1)$-form on $X_\CC$
\end{notation}

We fix, for future use, a choice of an algebraic closure $\Kbar$ of
$K$ and an embedding of $\Kbar$ into $\CC$.
% Also, we set $c_1(\L) := c_1(L)$ and
% $c_1(\overline {\mathcal L}) := c_1(\overline{L})$.

\subsection{Heights of cycles}
                  
% $\bullet$ {{\bf The height of cycles.}} 

Let $X$ be a smooth projective variety over $K$ of dimension
$N\geq 1$, and let $\mathcal X$ be a model for $X$ over $\OO_K$, i.e.,
an integral scheme projective and flat over Spec $\OO_K$ whose generic
fiber is $X$.

% For $p \geq 0$, let $\widehat{CH}^{p} (\mathcal {X})$ be the $p^{th}$
% \emph{arithmetic Chow group} of codimension $p$ arithmetic cycles
% on $\mathcal X$.
% Then we have the arithmetic intersection 
% pairing and the arithmetic degree map: 
% for each $p \geq 0$,
% $$
% \widehat{CH}^{p} (\mathcal {X}) \otimes \widehat{CH}^{1} (\mathcal {X})
% \longrightarrow \widehat{CH}^{p + 1} (\mathcal {X}),
% $$
% \noindent
% andp
% $$
% \widehat{\textup{deg}}: \widehat{CH}^{N + 1} (\mathcal {X}) 
% \longrightarrow \mathbb {R}.
% $$
% \noindent
% We identify an element in 
% $\widehat{CH}^{N + 1} (\mathcal {X})$ 
% and its degree via $\widehat{\deg}$. 

%\vspace{0.2cm}

% A hermitian metric is always assumed 
% to be smooth and invariant under complex conjugation.
% Let $\widehat {\textup{Pic}} \; (\mathcal {X})$ be the group of
% hermitian line bundles on $\mathcal X$ modulo isometric isomorphisms.
% Then we have the group homomorphism:
% $$
% \widehat{c_1} : \widehat {\textup{Pic}} \; (\mathcal {X}) \longrightarrow
% \widehat{CH}^{1} (\mathcal {X}).
% $$ 
% \noindent
% For a more thorough summary of the key properties of arithmetic Chow groups,
% see~\cite{Abb;97} or \cite{SUZ;97}.
% Proofs of the relevant facts can be found in ~\cite{GS;90}.

\vspace{0.2cm}

Let $\overline {\mathcal L}$ be a 
hermitian line bundle on $\mathcal X$. 
A hermitian metric is always assumed 
to be smooth and invariant under complex conjugation.
We assume furthermore that
$\mathcal L_K$ is ample,
and that the curvature form $c_1(\overline {\mathcal L})$ satisfies
$c_1(\overline {\mathcal L}) > 0$.  (See \cite{GS;90} for a discussion
of the curvature form associated to a hermitian line bundle).

\vspace{0.2cm}

By using arithmetic intersection theory, one defines the height of a cycle
in Arakelov geometry as follows (see e.g. \cite{SUZ;97}):

\begin{defn}
The \emph{height} of a nonzero effective cycle $Y$ (of pure
dimension) of $X$ with respect to 
$\overline {\mathcal L}$ is 
$$
h_{\overline {\mathcal L}} (Y) \ := \
{ {\widehat {c_1}(\overline {\mathcal L}|_{\overline {Y}})^
{\textup{dim} Y + 1}} \over
{ {(\textup{dim} Y + 1)} {c_1 ({\mathcal L}|_Y)^{\textup{dim} Y} }} },
$$
where $\overline {Y}$ is the (scheme-theoretic)
Zariski closure of $Y$ in $\mathcal {X}$.
\end{defn}

For a detailed overview of all the properties of curvature forms,
arithmetic Chern classes, and heights of arithmetic
cycles which we will need, see~\cite{Abb;97} or \cite{SUZ;97}.
Proofs of the relevant facts can be found in ~\cite{GS;90}.

\vspace{0.2cm}

% $\bullet$ {\bf {The cubical metric and the canonical height.}} 

\subsection{Canonical heights on abelian varieties}

\vspace{0.2cm}

Let $A/K$ be an abelian variety.
Using the choice of an embedding of $\Kbar$ into $\CC$, we view
$A(\Kbar)$ as a subset of $A(\CC)$.
Let $\mathcal A$ be a model for $A$ over Spec $\OO_K$.
Let $\overline {\mathcal {L}}$ be a hermitian line bundle on $\mathcal A$
such that 
$L := {\mathcal L}_K$ is symmetric and ample, and 
such that ${\mathcal L}$ is equipped with the \emph{cubical metric}
(see \cite{MB;85}).

Fix a nontrivial multiplication map on $A$ (e.g. multiplication by 2).
One can then construct 
from $(\mathcal A, \overline {\mathcal {L}})$ a sequence 
$({\mathcal A}_n, \overline{\mathcal {L}}_n)_{n \geq 1}$ of models of 
$(A, L)$, where each $\mathcal {L}_n$ 
is equipped with the cubical metric,
in such a way that the sequence
$$
h_{\overline {\mathcal L}_n} (Y) \ = \
{ {\widehat {c_1}(\overline {\mathcal L}_n|_{\overline {Y}_n})^
{\textup{dim} Y + 1}} \over
{ {(\textup{dim} Y + 1)} {c_1 ({L}|_Y)^{\textup{dim} Y} }}}
$$
converges (uniformly in $Y$) to a nonnegative real number $\hhat_L(Y)$. 
(Here $Y$ is a nonzero effective cycle
of pure dimension on $A$, and 
$\overline {Y}_n$ 
is the Zariski closure of $Y$ in $\mathcal {A}_n$.  
See ~\cite{SUZ;97} or~\cite{Zha;95;1} for details.)
% We denote the limit by $\hhat_L (Y)$. 
Note that $\hhat_L$ does not depend on the choice of 
$(\mathcal A, \mathcal L)$ or
a sequence of models 
$({\mathcal A}_n, \mathcal {L}_n)_{n \geq 1}$.
For $x \in A(\Kbar)$,   
$\hhat_L (x)$ is the N\'eron-Tate canonical height of $x$ with respect to $L$. 
% We also write $\overline L$ for $L$ with the cubical metric.

% For convenience, we often identify a closed subvariety (esp. point) of
% $A_{\overline {\mathbb Q}}$ 
% (or $A ({\overline {\mathbb Q}})$; 
% $A_{{\mathbb C}}$ or $A({{\mathbb C}})$)
% with its image in $A$ under the morphism
% $A_{\overline{\mathbb Q}} \rightarrow A$
% in case there is no danger of confusion
% or no specification makes any difference.
% Otherwise, we will specify what we mean.

\subsection{Statement of the main theorem}

We need several definitions in order to state our main result.  By a
\emph{variety} $X$ over a field $k$, we mean an integral separated
scheme of finite
type over $k$.  By a \emph{subvariety} of $X$, we mean an integral closed subscheme of $X$.

\begin{defn}
A \emph{torsion subvariety} of $A$ is
a translate of an abelian subvariety of $A$ by a torsion point.
% A closed subvariety of $A$ which is not
% torsion is called a \emph{nontorsion subvariety}.
\end{defn}

\begin{defn}
A sequence $(X_n)_{n \geq 1}$ of closed subvarieties of $A$ 
is \emph{small} if 
$\hhat_L(X_n) \rightarrow 0$ (as $n \to \infty$).
\end{defn}

% \begin{defn}
% Let $S$ be a closed subset of $A(\mathbb C)$,
% and let $A_0$ be a torsion subvariety of $A$ with $A_0(\mathbb C)
% \subseteq S$.  We say that $A_0$ is a \emph{maximal torsion
% subvariety} of $S$ if $A_0$ is
% the only torsion
% subvariety of $A$ that contains $A_0$ and has $A_0(\mathbb C)
% \subseteq S$.
% \end{defn}

\begin{defn}
A sequence $(X_n)_{n \geq 1}$ of closed subvarieties of $X$ is 
\emph{generic} if it has no subsequence contained in a
proper Zariski closed subset of $X$.
\end{defn}

\begin{defn}
A sequence $(X_n)_{n \geq 1}$ of closed subvarieties of $A$ is 
\emph{strict} if it has no subsequence contained in a
proper torsion subvariety of $A$.
\end{defn}

Note that the subvarieties $X_n$ are required to be defined over $\Kbar$, but not necessarily over $K$.

\vspace{0.2cm}

The following result is a generalization of
the Szpiro-Ullmo-Zhang/Ullmo/Zhang equidistribution theorem to
sequences of small subvarieties of an abelian variety. 

\vspace{0.2cm}

\begin{thm}[{\bf {Strict Equidistribution}}]~\label{thm;se}
Let $A/K$ be an abelian variety,
let $L$ be a symmetric ample line bundle on $A$, and
let $\overline L$ denote $L$ with the cubical metric.
Let $(X_n)_{n \geq 1}$ be a small
strict sequence of closed subvarieties of $A$. Then for every 
real-valued continuous
function $f$ on $A(\mathbb C)$, we have
$$
\int_{A ({\mathbb C})} f \; \mu_n
\longrightarrow
\int_{A ({\mathbb C})} f \; \mu,
$$
as $n \rightarrow \infty$, where 
setting
$d_n = \textup{dim} X_n$ and $g = \textup{dim} A$,
we have
$\mu_n = {\frac { {1} }{ { c_1({L|_{ X_{n}})^{d_n}} } }
{c_1 ({\overline L})^{d_n}} \delta_{X_n}}$ and
$\mu = { \frac {c_1 ({\overline L})^{g}}{c_1 (L)^{g}} }$.
%$\overline L$ is $L$ with the cubical metric, and
%$A ({\mathbb C}) 
%= \coprod_{ \sigma: k \hookrightarrow \mathbb {C} } A_{\sigma} (\mathbb {C})$
%if $A$ is defined over a number field 
%$k$ and $A_{\sigma}$ is its base extension to 
%$\mathbb {C}$ via $\sigma$. 
\end{thm}

{{\bf Remarks.}}
 
1. The first integral is the integral of $f$ against 
the restriction of ${c_1 ({\overline L})^{d_n}}/{\rm deg}_L (X_n)$ to $X_n(\CC)$.
The second integral is the integral of $f$ with respect to
the Haar measure $\mu$ on $A(\mathbb {C})$, normalized to have total
mass 1.

2. If $X_n = x_n$ is a point, i.e., if $d_n = 0$, note that
$$
\int_{A ({\mathbb C})} f(x) \mu_n \ = \
{ \frac{1}{\# O(x_n)} } \sum_{x \in O(x_n)} f(x),
$$
\noindent
where $O(x_n)$ is the orbit of $x_n$ under the action of 
$\Gal(\Kbar / K)$.

% There is some abuse of notation here:
% the $x_n$ in the definition of $\mu_n$ stands for a closed point of
% $A$ (i.e., a Galois conjugacy class of points of $A(\Kbar)$), 
% while the $x_n$ in $O(x_n)$ stands for one of the corresponding points in $A(\mathbb C)$.

3. For notational convenience, we write
$$
\mu_n \;
\stackrel{w}{\longrightarrow} \;  \mu \;\; \mathrm{as} \;\; 
n \rightarrow \infty,
$$  
and say the sequence $(\mu_n)_{n \geq 1}$ of measures 
\emph{weakly converges} to $\mu$, if 
$$
\int_{A({\mathbb C})} f \mu_n \to \int_{A({\mathbb C})} f \mu
$$
for every continuous function $f : A({\mathbb C}) \to {\mathbb R}$.
In this case, we say that the
$X_n$'s are \emph{equidistributed} with respect to $\mu$. 

4. To get a feeling for what Theorem~\ref{thm;se} says, consider the following
simple example.
Let $E$ be an elliptic curve defined over $\mathbb Q$ and let 
$A = E \times E$. For each $n \geq 1$, let 
$E_n \subset A$ be the graph of the multiplication-by-$n$
map on $E$.  Then each $E_n$ is a torsion subvariety 
of $A$ defined over $\mathbb Q$ (in fact, $E_n$ is $\mathbb
Q$-isogenous to $E$).  
% Thus there are infinitely many torsion subvarieties of $A$ which are defined 
% over ${\mathbb Q}$,
% cf. Corollary 2,~\cite{Zha;98}. 

It is easy to see that ${\rm deg}(E_n) \to \infty$ as $n\to\infty$,
and that $\bigcup_{n\geq 1} E_n$ is Zariski dense in $A$.
Theorem~\ref{thm;se} says something stronger than this, namely that as
$n\to\infty$, the normalized Haar measure on $E_n$ approximates the
normalized Haar measure on $A$ arbitrarily closely.

5. For related equidistribution results, see
Theorem 1.1 of~\cite{Bil;97},
Theorem 4.1 of~\cite{SUZ;97},
Theorem 2.3 of~\cite{Ull;98}, and
Theorem 1.1 of~\cite{Zha;98}.
% In addition, the authors recently learned that Pascal Autissier has
% independently obtained a proof of Theorem~\ref{thm;ge} below.
In addition, 
Pascal Autissier has recently obtained a proof of Theorem~\ref{thm;ge}
of the present paper independently of the authors.

% Similar remarks will also apply to
% a similar situation in Generic Equidistribution 
% Theorem~\ref{thm;ge} below
% without mentioning it explicitly.

\vspace{0.2cm}

% The following corollary is motivated by some remarks in~\cite{Zha;98}. 
% Let $S$ be a proper closed subset of $A (\mathbb C)$, and
% let $G := \Gal ({\overline {\mathbb Q}}/{\mathbb Q})$.
% Applying Theorem~\ref{thm;se} to a nonzero
% real-valued continuous function $f \geq 0$ on $A(\mathbb C)$ whose
% support is disjoint from $S$, 
% one easily deduces the following result:
%
% \vspace{0.2cm}
%
% \begin{cor}~\label{cor;epsilon}
% Let $S$ be a proper closed subset of $A (\mathbb C)$ 
% in the complex topology. Then:
%
% (i) \hspace{0.01cm} There are only finitely many maximal torsion
% subvarieties $A_0$ of $S$ with $A_0^\sigma \subseteq S$ for all $\sigma
% \in G$.  
%
% (ii) Let $A'$ be the union of all maximal torsion
% subvarieties $A_0$ of $S$ with $A_0^\sigma \subseteq S$ for all $\sigma
% \in G$, and let $S^{\prime}$ be the complement of $A'$ in $S$.
% Then there is a real number $\epsilon > 0$ such that $\hhat_L (Y) > \epsilon$
% for every closed subvariety $Y$ of $A$ with $Y^\sigma \subseteq
% S^\prime$ for all $\sigma \in G$.
% 
% \end{cor}

\section{Generic Equidistribution}~\label{sec:ge}

\vspace{0.2cm}

Let $X$ be a closed subvariety of
dimension $N \geq 1$ of $A$. 
The following result is a special case of Zhang's 
``Theorem of the successive minima''(see~\cite{Zha;98} for details):
\begin{thm}
\label{SuccesiveMinimaTheorem}
Define 
$$
\lambda_1 (X) := \sup_{Z} \inf_{x \in X - Z} \hhat_L (x),
$$
\noindent
where $Z$ runs over the set of all proper closed subsets of $X$.
% For $i = 1,$ 2, $\dots$, $N + 1$, let
% $$
% \lambda_i (X) := \sup_{Z} \inf_{x \in X - Z} \hhat_L (x),
% $$
% \noindent
% where $Z$ runs over the set of codimension $i$ closed subsets of
% $X$.
Then
$$
\lambda_1 (X) \ \geq \ \hhat_L(X) 
\ \geq \ 
{ \frac{1}{N + 1} } \lambda_1 (X).
% { {1} \over {N + 1} } \sum_{i=1}^{N+1} \lambda_i (X).
$$
\end{thm}
\noindent 

\begin{defn}
Let $\overline{\L}$ be a hermitian line bundle on $\A$.  
If $f$ is a real-valued $C^{\infty}$-function on $A({\mathbb C})$, define
$$
\overline {\mathcal L} (f) \ := \
\overline {\mathcal L} \otimes ({\mathcal O}_{\mathcal A}, e^{- f})
$$
to be the tensor product of $\overline{\L}$ with the trivial bundle,
endowed with the metric given by $\| 1 \|(P) = e^{-f(P)}$.
\end{defn}

\vspace{0.2cm}

\begin{thm}[{\bf {Generic Equidistribution}}]~\label{thm;ge}
Let $A/K$ be an abelian variety, and
let $L$ be a symmetric ample line bundle on $A$.
Let $(X_n)_{n \geq 1}$ be a small
generic sequence of closed subvarieties of $A$. Then, for every 
real-valued continuous
function $f$ on $A(\mathbb C)$, we have
$$
\int_{A ({\mathbb C})} f(x) \mu_n 
\longrightarrow
\int_{A ({\mathbb C})} f(x) \mu
$$
as $n \rightarrow \infty$, where 
$d_n = \textup{dim} X_n$,
$\mu_n = 
{\frac { {1} }{ { c_1({L|_{ X_{n}})^{d_n}} } }
{c_1 ({\overline L})^{d_n}} \delta_{X_n}}$,
$g = \textup {dim} A$,
$\mu = { \frac{c_1 ({\overline L})^{g}}{c_1 (L)^{g}} }$, 
and $\overline L$ is $L$ with the cubical metric. 
\end{thm}

\emph{Proof}.
% Without loss of generality, we may assume that the $X_n$'s are distinct.
% By rearranging the $X_n$'s if necessary, 
% we may furthermore assume that 
% $X_j \not\supset X_i$ if $j < i$. 
Enumerate the countably many subvarieties $( Z_n )_{n \geq 1}$ of $A$ defined
over $\overline{K}$.  Since $(X_n)_{n \geq 1}$ is generic, we may assume, 
%(after passing to a subsequence if necessary) 
without loss of generality, 
$X_n \not\subset Z_1 \cup \cdots \cup Z_n$.  
By the definition of $\lambda_1(X_n)$, we can find (for each $n\geq
1$) an infinite sequence 
%$(x_{n,m}) \in X_n$ 
$(x_{n,m})_{m \geq 1}$ in $X_n$ 
such that:

\vspace{0.1cm}

(i) \hspace{0.08cm} for each $m \geq 1$,
$x_{n,m} \notin \bigcup_{1 \leq i \leq n} Z_i $; 

(ii) $|\hhat_L(x_{n,m}) - \lambda_1 (X_n)| <
 \frac{1}{n}$ for all $m \geq 1$; \ \ \and

(iii) for each $n\geq 1$, $\lim_{m\to\infty} \hhat_L(x_{n,m}) =
 \lambda_1 (X_n)$.

\vspace{0.1cm}

By choosing a bijection between $\NN^2$ and $\NN$, we may consider the
doubly-indexed sequence $(x_{n,m})$ as a sequence indexed by the natural numbers.
Property (i) guarantees that the resulting sequence $(x_{n,m})$ is generic.
Furthermore, since $\hhat_L(X_n) \to 0$ by assumption, it follows from the theorem of the
successive minima that $\lambda_1(X_n) \rightarrow 0$
as $n \rightarrow \infty$.
Using this observation, properties (ii) and (iii) easily imply that $\lim_{n, m \to \infty}
\hhat_L(x_{n,m}) = 0$, i.e., that the sequence $(x_{n,m})$ is small.

\vspace{0.1cm}

Define
$$
\alpha_{n,m} :=
{ \frac{1}{\# O(x_{n,m})} } \sum_{x \in O(x_{n,m})} f(x),
$$
where $O(x_{n,m})$ is the orbit of $x_{n,m}$ under the action of
$\Gal(\Kbar / K)$.
By choosing a
subsequence of $(x_{n,m})_{m \geq 1}$ if necessary, we may assume that
$\lim_{m \rightarrow \infty} \alpha_{n,m}$ exists for all $n \geq 1$.
Note that every subsequence of a small (resp. generic) sequence is
small (resp. generic).

Approximating $f$ by $C^{\infty}$-functions if necessary, we may
assume that $f$ is a $C^{\infty}$-function. 
Let $\lambda > 0$ be a real number. Note that 
$c_1(\overline{\mathcal L}_l (\lambda f)) > 0$ if $\lambda >0$ is small 
enough.
We then note, for $l \geq 1$, that
\begin{eqnarray}
h_{ \overline{\mathcal L}_l (\lambda f) } (x_{n,m}) 
\ & = & \ 
h_{ \overline{\mathcal L}_l } (x_{n,m}) + 
\lambda \alpha_{n,m}; \  \  \;\;
\textup{and} \nonumber\\
\liminf_{m \rightarrow \infty} 
h_{ \overline{\mathcal L}_l (\lambda f) } (x_{n,m}) 
\ & \geq & \ 
h_{ \overline{\mathcal L}_l (\lambda f) } (X_n) \ \ \ \ \ \ \ \ \ \ \ \ \ \;
\textup{\cite{SUZ;97}, Proposition 2.1}
\nonumber\\
\ & = & \
h_{ \overline{\mathcal L}_l } (X_n) +   
\lambda \int_{A (\mathbb C)} f(x) \mu_n + O(\lambda^2), \nonumber 
\end{eqnarray}
\noindent
where the last equality follows from \cite[Proof of Proposition 2.9]{Abb;97}.
Here the $O$-constant is independent of $l$ and $n$, and $\lambda > 0$
is sufficiently small.

Fix $n\geq 1$ and $\varepsilon > 0$. 
Then for $m$ sufficiently large, we have:
$$
h_{ \overline{\mathcal L}_l } (x_{n,m}) + 
\lambda \alpha_{n,m}
\  \geq  \ 
h_{ \overline{\mathcal L}_l } (X_n) +   
\lambda \int_{A (\mathbb C)} f(x) \mu_n + O(\lambda^2) - \varepsilon.
$$

\noindent Letting $l \rightarrow \infty$, we have:
$$
\hhat_L (x_{n,m}) - \hhat_L (X_n) +
\lambda \alpha_{n,m} \ \geq \ 
\lambda \int_{A (\mathbb C)} f(x) \mu_n +
O(\lambda^2) - \varepsilon.
$$

\noindent Now let $m \rightarrow \infty$, and we obtain (since
$\epsilon > 0$ is arbitrary):
\begin{eqnarray}~\label{eqnarray;m}
\lambda_1 (X_n) - \hhat_L (X_n) +
\lambda \lim_{m \rightarrow \infty} \alpha_{n,m}
% \lim_{m \rightarrow \infty} \alpha_{n,m} 
\ \geq \ 
\lambda \int_{A (\mathbb C)} f(x) \mu_n +
O(\lambda^2). 
\end{eqnarray}

% Note any subsequence $(x_{n_j, m_i})_{j, i \geq 1}$ 
% of $(x_{n, m})_{n, m \geq 1}$ 
% is a small generic sequence in $A$.
% Recall it means, in particular, that 
% $\lim_{j, i \rightarrow \infty} \hhat_L(x_{n_j,m_i}) = 0$.

On the other hand, the Szpiro-Ullmo-Zhang/Ullmo/Zhang equidistribution theorem
~(\cite{SUZ;97},~\cite{Ull;98}, and~\cite{Zha;98}), applied to 
the small generic sequence $(x_{n,m})$, implies that
$\lim_{n,m \to \infty} \alpha_{n,m}$ exists, and that
\begin{eqnarray}~\label{eqnarray;suz}
\lim_{n, m \rightarrow \infty} \alpha_{n,m}
\ = \
\int_{A(\mathbb C)} f(x) \mu. 
\end{eqnarray}

Taking $\limsup_{n \rightarrow \infty}$ in~(\ref{eqnarray;m}),
we have:
\begin{eqnarray}
\lambda \lim_{n,m \rightarrow \infty}
\alpha_{n,m} 
\ \geq \ 
\lambda \limsup_{n \rightarrow \infty} 
\int_{A (\mathbb C)} f(x) \mu_n +
O(\lambda^2). \nonumber
\end{eqnarray}

\noindent 
Now divide both sides by $\lambda > 0$, and let $\lambda \rightarrow
0^+$. We obtain:
\begin{eqnarray}~\label{eqnarray;limsup}
\lim_{n,m \rightarrow \infty} \alpha_{n,m} 
\ \geq \ 
\limsup_{n \rightarrow \infty} 
\int_{A (\mathbb C)} f(x) \mu_n. 
\end{eqnarray}

\noindent Replacing $f$ by $-f$, we see that:
\begin{eqnarray}~\label{eqnarray;liminf}
\lim_{n,m \rightarrow \infty} \alpha_{n,m} 
\ \leq \ 
\liminf_{n \rightarrow \infty} 
\int_{A (\mathbb C)} f(x) \mu_n. 
\end{eqnarray}

It then follows from~(\ref{eqnarray;limsup}) and~(\ref{eqnarray;liminf})
that 
$\lim_{n \rightarrow \infty} \int_{A (\mathbb C)} f(x) \mu_n$ exists and
that
$$
\lim_{n \rightarrow \infty} 
\int_{A (\mathbb C)} f(x) \mu_n 
\ = \ 
\lim_{n,m \rightarrow \infty} \alpha_{n,m}.
$$

\noindent
We conclude from~(\ref{eqnarray;suz}) that
$$
\lim_{n \rightarrow \infty} 
\int_{A (\mathbb C)} f(x) \mu_n \ = \
\int_{A(\mathbb C)} f(x) \mu, 
$$ 
\noindent as desired. $\hfill \square$

\vspace{0.2cm}

\section{Strict Equidistribution}
~\label{sec:se}

\vspace{0.2cm}

The following result is a consequence of two results of Zhang: the generalized Bogomolov
conjecture (see \cite{Zha;98}) and the theorem of the successive
minima.  The proof is similar to the proof of Theorem~\ref{thm;ge}.

\begin{thm}
% [{\bf {Zhang~\cite{Zha;98}: 
% Generalized Bogomolov Conjecture}}]~
\label{thm;bogo} 
Let $X$ be a nontorsion
subvariety of $A$. Then there is an $\epsilon > 0$ such that the set
$$
\bigcup \big \{ Y : Y \; \textup{is a closed subvariety of} \; X \; 
\textup{such that} 
\; \hhat_L (Y) 
\leq \epsilon \big \}
$$
is not Zariski dense in $X$.
\end{thm}

\emph{Proof}. 
% (See also the beginning of the proof of Theorem~\ref{thm;ge}.)
Suppose, for the sake of contradiction, that $(Y_n)_{n \geq 1}$ is a
sequence of distinct closed subvarieties of $X$ which is
small (i.e., $\hhat_L (Y_n) \rightarrow 0$) and  generic in
$X$ (i.e., no subsequence is contained in a proper Zariski closed subset of
$X$).
Then, proceeding as in the proof of Theorem~\ref{thm;ge}, we can
construct an infinite sequence $(y_k)_{k \geq 1}$ of points in $X$ such that 
$\{ y_k \in X: k \geq 1 \}$ is Zariski dense in $X$
and $\hhat_L(y_k) \rightarrow 0$.
But Corollary 3 of~\cite{Zha;98} then implies  
that $X$ is a torsion subvariety of $A$, a contradiction. $\hfill \square$

\vspace{.2cm}

Now we are ready to prove Theorem~\ref{thm;se} (Strict Equidistribution 
Theorem).

\vspace{0.2cm}

\emph{Proof} (of Theorem~\ref{thm;se}). By Theorem~\ref{thm;ge},
it suffices to show that the small and strict sequence 
$(X_n)_{n \geq 1}$ is generic.  
Let $X'$ be the Zariski closure of $\bigcup_k X_{n_k}$ for any
subsequence $(X_{n_k})_{k \geq 1}$ of $(X_n)_{n \geq 1}$.  
By Theorem~\ref{thm;bogo}, 
$X$ must be a torsion subvariety of $A$.  Since $(X_n)_{n\geq 1}$ is
strict, it follows that $X' = A$, so that $(X_n)_{n\geq 1}$ is generic
as desired.
$\hfill \square$

\vspace{0.2cm}

% {{\bf Remark.}}
% A similar proof gives an analogous equidistribution result 
% in the more general context of an \emph{adelic metric},
% cf.~\cite{Zha;98}. 

% 1. One can also prove Theorem~\ref{thm;se} directly
% for small strict sequences in a way very analogous to the proof of
% Theorem~\ref{thm;ge}, yet still together with the use of
% Theorem~\ref{thm;bogo}. There is no big difference at all.   

%%%%%%%%%%%%%%%%%%%%%%%%%%%%%%%%%%%%%%%%%%%%%%%%%
%               Bibliography                    %
%%%%%%%%%%%%%%%%%%%%%%%%%%%%%%%%%%%%%%%%%%%%%%%%%

\bibliographystyle{amsplain}

\begin{thebibliography}{10}
\addcontentsline{toc}{chapter}{Bibliography}

\bibitem{Abb;97}
A.~Abbes, \emph{Hauteurs et discr\'etude 
(d'apr\`es L. Szpiro, E. Ullmo et S. Zhang)}, 
(French) [Heights and discreteness (after L. Szpiro, E. Ullmo and S. Zhang)] 
S\'eminaire Bourbaki, Vol. 1996/97. Ast\'erisque No. {{\bf 245}} (1997), 
Exp. No. {{\bf 825}}, 4, 141--166.

\bibitem{Bil;97}
Y.~Bilu, \emph{Limit distribution of small points on algebraic tori}, Duke Math. J. \textbf{89} (1997), no. 3, 465--476. 

\bibitem{GS;90}
H.~Gillet, and C.~Soul\'e, \emph{Arithmetic intersection theory}, Inst. 
Hautes \'Etudes Sci. Publ. Math. No. \textbf{72} (1990), 93--174.

\bibitem{MB;85} 
L.~Moret-Bailly, \emph{M{\'e}triques permises} (French) [Admissible metrics] S\'eminaire sur les pinceaux de
  courbes elliptiques, L.~Szpiro (ed.), Ast\'erisque No. {{\bf 127}}
  (1985), 29--87.

\bibitem{SUZ;97}
L.~Szpiro, E.~Ullmo, and S.~Zhang, \emph{\'Equir\'epartition des petits points} (French), [Uniform distribution of small points] Invent. Math. \textbf{127} (1997), no. 2, 337--347. 

\bibitem{Ull;98}
E.~Ullmo, \emph{Positivit\'e et discr\'etion des points alg\'ebriques des courbes} (French), [Positivity and discreteness of algebraic points of curves] Ann. of Math. (2) \textbf{147} (1998), no. 1, 167--179. 

\bibitem{Zha;95}
S.~Zhang, \emph{Positive line bundles on arithmetic varieties}, J. Amer. Math. Soc. \textbf{8} (1995), no. 1, 187--221.

\bibitem{Zha;95;1}
\bysame, \emph{Small points and adelic metrics}, J. Algebraic Geom. 
{{\bf 4}} (1995), no. 2, 281--300.

\bibitem{Zha;98}
\bysame, \emph{Equidistribution of small points on abelian varieties},
Ann. of Math. (2) \textbf{147} (1998), no. 1, 159--165.




%\vspace{4cm}

\end{thebibliography}

\providecommand{\bysame}{\leavevmode\hbox to3em{\hrulefill}\thinspace}

\end{document}